\documentclass[11pt]{amsart}
\usepackage{graphicx, amscd, amsfonts, amsthm}
\title{Homotopy types of homeomorphism groups of noncompact 2-manifolds}
\author{Tatsuhiko Yagasaki}
\address{Department of Mathematics, Kyoto Institute of Technology, Matsugasaki, Sakyoku, Kyoto 606, Japan}
\email{yagasaki@ipc.kit.ac.jp}

\newtheorem{theorem}{Theorem}[section]
\newtheorem{proposition}{Proposition}[section] 
\newtheorem{corollary}{Corollary}[section] 
\newtheorem{lemma}{Lemma}[section]
\newtheorem*{claim}{Claim}

\theoremstyle{definition}

\newtheorem{fact}{Fact}[section]

\newtheorem*{Proof of Theorem 3.1}{Proof of Theorem 3.1}
\newtheorem*{Proof of Corollary 3.1}{Proof of Corollary 3.1}
\newtheorem*{Proof of Proposition 4.1.(1)}{Proof of Proposition 4.1.(1)}
\newtheorem*{Proof of Proposition 4.1.(2)}{Proof of Proposition 4.1.(2)}
\newtheorem*{Proof of Corollary 1.1}{Proof of Corollary 1.1}

\def \cal {\mathcal}
\def \phi {\varphi}

\begin{document}
\baselineskip 6 mm

\thispagestyle{empty}

\maketitle
\begin{abstract}
Suppose $M$ is a noncompact connected PL 2-manifold and let ${\cal H}(M)_0$ denote the identity component of the homeomorphism group of $M$ with the compact-open topology. In this paper we classify the homotopy type of ${\cal H}(M)_0$ by showing that ${\cal H}(M)_0$ has the homotopy type of the circle if $M$ is the plane, an open or half open annulus, or the punctured projective plane. In all other cases we show that ${\cal H}(M)_0$ is homotopically trivial.
\end{abstract}

%Section 1

\section{Introduction}

Hamstrom \cite{Ha} classified the homotopy types of the identity components of homeomorphism groups of compact 2-manifolds $M$. In this paper we treat the case where $M$ is noncompact. Suppose $M$ is a PL 2-manifold and $X$ is a compact subpolyhedron of $M$. We denote by ${\cal H}_X(M)$ the group of homeomorphisms $h$ of $M$ onto itself with $h|_X =id$, equipped with the compact-open topology, and by ${\cal H}(M)_0$ the identity component of ${\cal H}(M)$. Let ${\Bbb R}^2$ denote the plane, ${\Bbb S}^1$ the unit circle and ${\Bbb P}^2$ the projective plane. The following is the main result of this paper.

\begin{theorem}
Suppose $M$ is a noncompact connected (separable) PL $2$-manifold and $X$ is a compact subpolyhedron of $M$. Then \\
{\rm (i)} ${\cal H}_X(M)_0 \simeq {\Bbb S}^1$ if $(M, X) \cong ({\Bbb R}^2, \emptyset)$, $({\Bbb R}^2, 1pt)$, $({\Bbb S}^1 \times {\Bbb R}^1, \emptyset)$, $({\Bbb S}^1 \times [0, 1), \emptyset)$ or $({\Bbb P}^2 \setminus 1pt, \emptyset)$, \\
{\rm (ii)} ${\cal H}_X(M)_0 \simeq \ast$ in all other cases.  
\end{theorem}

\begin{corollary}
If $M$ is a connected (separable) $2$-manifold and $X$ is a compact subpolyhedron of $M$ with respect to some triangulation of $M$, then ${\cal H}_X(M)_0$ is an $\ell_2$-manifold. 
\end{corollary}

In \cite{Ya} we obtained a natural principal bundle connecting the homeomorphism group and the embeddding space (cf.\,Section 2). In this paper we will seek a condition under which the fiber of this bundle is connected (Section 3). The contractibility and the ANR property of ${\cal H}_X(M)_0$ in the compact case will then imply the similar properties of embedding spaces and in turn the corresponding properties of ${\cal H}_X(M)_0$ in the noncompact case. Corollary 1.1 follows immediately from the characterization of $\ell_2$-manifolds and this enables us to determine the topological type itself of ${\cal H}_X(M)_0$ by the homotopy invariance of infinite-dimensional manifolds.

In a succeeding paper we will investigate the subgroups of ${\cal H}_X(M)_0$ consisting of PL and Lipschitz homeomorphisms  from the viewpoints of infinite-dimensional topological manifolds. 

\section{Preliminaries}

Throughout the paper we follow the following conventions: Spaces are assumed to be separable and metrizable, and maps are always continuous. When $A$ is a subset of a space $X$, the notations ${\rm Fr}_X A$, ${\rm cl}_X A$ and ${\rm Int}_X A$ denote the frontier, closure and interior of $A$ relative to $X$ (i.e., ${\rm Int}_X A = \{ x \in A \, | \, A$ contains a neighborhood of $x$ in $X \}$ and ${\rm Fr}_X A = {\rm cl}_X A \setminus {\rm Int}_X A$). On the other hand, when $M$ is a manifold, the notations $\partial M$ and ${\rm Int}\,M$ denote the boundary and interior of $M$ as a manifold. 
When $N$ is a 2-submanifold of a 2-manifold $M$, we always assume that $N$ is a closed subset of $M$ and ${\rm Fr} N = {\rm Fr}_M N$ is a 1-manifold transversal to $\partial M$. Therefore we have ${\rm Int}\,N = {\rm Int}_M N \cap {\rm Int}\,M$ and ${\rm Fr}_M N \subset \partial N$. 
% We simply write as ${\rm Fr}\,N$ when the meaning is clear from the context. 
A metrizable space $X$ is called an ANR (absolute neighborhood retract) if any map $f : B \to X$ from a closed subset of a metrizable space $Y$ has an extension to a neighborhood $U$ of $B$. If we can always take $U = Y$, then $X$ is called an AR (absolute retract). ANRs are locally contractible and ARs are exactly contractible ANRs (cf.\,\cite{Hu}). Finally $\ell_2$ denotes the separable Hilbert space $\{ (x_n) \in {\Bbb R}^\infty : \sum_n x_n^2 < \infty \}$. 

In \cite{Ya} we investigated some extension property of embeddings of a compact 2-polyhedron into a 2-manifold, based upon the conformal mapping theorem. The result is summarized as follows:  
Suppose $M$ is a PL 2-manifold and $K \subset X$ are compact subpolyhedra of $M$. Let ${\cal E}_K(X, M)$ denote the space of embeddings $f : X \hookrightarrow M$ with $f|_K = id$, equipped with the compact-open topology. We consider the subspace of proper embeddings ${\cal E}_K(X, M)^{\ast} = \{f \in {\cal E}_K(X,M) \, : \, f(X \cap \partial M) \subset \partial M, f(X \cap {\rm Int} \,M) \subset {\rm Int} \, M \}$. Let ${\cal E}_K(X, M)^{\ast}_0$ denote the connected component of the inclusion $i_X : X \subset M$ in ${\cal E}_K(X, M)^{\ast}$.  

\begin{theorem}
For every $f \in {\mathcal E}_K(X, M)^{\ast}$ and every neighborhood $U$ of $f(X)$ in $M$, there exists a neighborhood ${\mathcal U}$ of $f$ in ${\mathcal E}_K(X, M)^{\ast}$ and a map $\phi : {\mathcal U} \to {\mathcal H}_{K \cup (M \setminus U)}(M)_0$ such that $\phi(g)f = g$ for each $g \in {\mathcal U}$ and $\phi(f) = id_M$.
\end{theorem}

\begin{corollary}
For any open neighborhood $U$ of $X$ in $M$, the restriction map $\pi : {\cal H}_{K \cup (M \setminus U)}(M)_0 \to {\cal E}_K(X, U)_0^{\ast}$, $\pi(f) = f|_X$, is a principal bundle with fiber ${\cal G} \equiv {\cal H}_{K \cup (M \setminus U)}(M)_0 \cap {\mathcal H}_X(M)$, where the group ${\mathcal G}$ acts on ${\cal H}_{K \cup (M \setminus U)}(M)_0$ by right composition. 
\end{corollary}

\begin{proposition}
${\cal E}_K(X, M)$ and ${\cal E}_K(X, M)^{\ast}$ are ANRs.
\end{proposition}

Next we recall some fundamental facts on homeomorphism groups of compact 2-manifolds. 

\begin{fact} 
If $N$ is a compact PL $2$-manifold and $Y$ is a compact subpolyhedron of $N$, then ${\cal H}_Y(N)$ is an ANR. (\cite{Ja}, \cite{LM}, cf.\,\cite[Lemma 3.2]{Ya}) 
\end{fact} 

\begin{lemma} {\rm (}\cite{Ha}, \cite[\S3]{Sc}{\rm )} Suppose $N$ is a compact connected PL $2$-manifold and $Y$ is a compact subpolyhedron of $N$. \\
{\rm (i)} If $(N, Y) \not\cong ({\Bbb D}^2, \emptyset)$, $({\Bbb D}^2, 0)$, $({\Bbb S}^1 \times [0,1], \emptyset)$, $({\Bbb M}, \emptyset)$, $({\Bbb S}^2, \emptyset)$, $({\Bbb S}^2, 1pt)$, $({\Bbb S}^2, 2pts)$, $({\Bbb T}^2, \emptyset)$, $({\Bbb K}^2, \emptyset)$, $({\Bbb P}^2, \emptyset)$, $({\Bbb P}^2, 1pt)$, then ${\cal H}_Y(N)_0 \simeq \ast$. \\
{\rm (ii)} If $A$ is a nonempty compact subset of $\partial \, N$, then ${\cal H}_{Y \cup A}(N)_0 \simeq \ast$. \\
{\rm (iii)} If $(N, Y) \cong ({\Bbb D}^2, \emptyset)$, $({\Bbb D}^2, 0)$, $({\Bbb S}^1 \times [0,1], \emptyset)$, $({\Bbb M}, \emptyset)$, $({\Bbb S}^2, 1pt)$, $({\Bbb S}^2, 2pts)$, $({\Bbb P}^2, 1pt)$, $({\Bbb K}^2, \emptyset)$, then ${\cal H}_Y(N)_0 \simeq {\Bbb S}^1$. 
\end{lemma}

\begin{proof}
In \cite{Sc} the PL-homeomorphism groups of compact 2-manifolds was studied in the context of semisimplicial complex. However, using Corollary 2.1 and the results in \cite{Ha}, we can apply the arguments and results in \cite[\S3]{Sc} to our setting. 

(i) Let $L$ be a small regular neighborhood of the union $Y_1$ of the nondegenerate components of $Y$ and let $Y_0 = Y \setminus Y_1$. Since ${\cal H}_Y(N)_0$ deforms into ${\cal H}_{L \cup Y_0}(N)_0 \cong {\cal H}_{{\rm Fr} \, L \cup Y_0}(cl \, (N \setminus L))_0$, we may assume that $Y_1 \subset \partial N$. This case follows from \cite{Ha} and \cite[\S3]{Sc}. 

(ii) Let $\partial_+ N$ denote the union of the components of $\partial \, N$ which meet $A$. Then ${\cal H}_{Y \cup A}(N)_0$ strongly deformation retracts onto ${\cal H}_{Y \cup \partial_+ N}(N)_0$, and the latter is contractible by the case (i).  
\end{proof}

\section{Relative isotopies on 2-manifolds} 

In Corollary 2.1 we have a principal bundle with a fiber ${\mathcal G} \equiv {\cal H}_X(M) \cap {\cal H}_K(M)_0$. In this section we will seek a sufficient condition which implies ${\mathcal G} = {\cal H}_X(M)_0$. Suppose $M$ is a 2-manifold and $N$ is a 2-submanifold of $M$. In \cite{Ep} it is shown that (i) two homotopic essential simple closed curves in ${\rm Int}M$ and two proper arcs homotopic rel ends in $M$ are ambient isotopic rel $\partial M$, (ii) every homeomorphism $h : M \to M$ homotopic to $id_M$ is ambient isotopic to $id_M$. Using these results or arguments we will show that if, in addition, $h|_N = id_N$ then $h$ is isotopic to $id_M$ rel $N$ under some restrictions on disks, annuli and M\"obius bands components (i.e. the pieces which admit global rotations). We denote the M\"obius band, the torus and the Klein bottle by ${\Bbb M}$, ${\Bbb T}^2$ and ${\Bbb K}^2$ respectively. The symbol $\#X$ denotes the number of elements (or cardinal) of a set $X$.

\begin{theorem}
Suppose $M$ is a connected $2$-manifold, $N$ is a compact $2$-submanifold of $M$ and $X$ is a subset of $N$ such that \\ {\rm (i)} $M \neq {\Bbb T}^2$, ${\Bbb P}^2$, ${\Bbb K}^2$ or $X \neq \emptyset$. \\
{\rm (ii)}{\rm (a)} if $H$ is a disk component of $N$, then $\# (H \cap X) \geq 2$, \\
{\rm (b)} if $H$ is an annulus or M\"obius band component of $N$, then $H \cap X \neq \emptyset$, \\
{\rm (iii)} {\rm (a)} if $L$ is a disk component of $cl(M \setminus N)$, then ${\rm Fr}L$ is a disjoint union of arcs or $\# (L \cap X) \geq 2$, \\ {\rm (b)} if $L$ is a M\"obius band component of $cl(M \setminus N)$, then ${\rm Fr}L$ is a disjoint union of arcs or $L \cap X \neq \emptyset$. \\ 
If $h_t : M \to M$ is an isotopy rel $X$ such that $h_0|_N = h_1|_N$, then there exists an isotopy $h'_t : M \to M$ rel $N$ such that $h'_0 = h_0$, $h'_1 = h_1$ and $h'_t = h_t \, (0 \leq t \leq 1)$ on $M \setminus K$ for some compact subset $K$ of $M$. 
\end{theorem}

\begin{corollary}
Under the same condition as in Theorem 3.1, we have ${\cal H}_N(M) \cap {\cal H}_X(M)_0 = {\cal H}_N(M)_0$. 
\end{corollary}

First we explain the meaning of the conditions (ii) and (iii) in Theorem 3.1. 
Suppose $h \in {\cal H}_N(M)$ and $h$ is isotopic to $id_M$ rel $X$. In order that $h$ is isotopic to $id_M$ rel $N$, it is necessary that $h$ does not Dehn twist along the boundary circle of any disk, M\"obius band and annulus component of $N$. 
This is ensured by the condition (ii) in Theorem 3.1 (Figure 1.$a$: $N = N_1 \cup N_2$, $cl(M \setminus N) = L$ and $X = \{ a, b\}$). 
However, this is not sufficient because a union of some components of $N$ and $cl(M \setminus N)$ may form a disk, a M\"obius band or an annulus. 
The condition (iii) in Theorem 3.1 is imposed to prevent Dehn twists around these pieces (Figure 1.$b$: $N = N_1 \cup N_2$, $cl(M \setminus N) = L \cup L_1 \cup L_2$ and $X = \{ a_1, a_2, b_1, b_2\}$). 
This condition is too strong (we can replace $X$ by $Y = \{ a_1, b_1 \}$), but it is simple and sufficient for our purpose.

\begin{figure}[h]
\caption{$h$ can not Dehn twist $L$.}
\[ \fbox{Figure 1.a} \hspace{2cm} \fbox{Figure 1.b} \]
\end{figure}

We proceed to the verification of Theorem 3.1. We need some preliminary lemmas. Throughout this section we assume that $M$ is a connected $2$-manifold and $N$ is a $2$-submanifold of $M$. When $G$ is a group and $S \subset G$, $\langle S \rangle$ denotes the subgroup of $G$ generated by $S$. 

We will use the fallowing facts from \cite{Ep}.

\begin{fact}
(0) (\cite[Theorem 1.7]{Ep}) If a simple closed curve $C$ in $M$ is null-homotopic, then it bounds a disk.  \\
(1) (\cite[Theorem 3.1]{Ep}) Suppose $\alpha$ and $\beta$ are proper arcs in $M$. If they are homotopic relative to end points, then they are ambient isotopic relative to $\partial M$.  \\
(2) (\cite[Theorem 4.2]{Ep}) Let $C$ be a simple closed curve in $M$, which does not bound a disk or a M\"obius band. Let $\alpha \in \pi_1(M,\ast)$ be represented by a single circuit of $C$ and let $\alpha = \beta^k$, $k \geq 0$. Then $\alpha = \beta$.  \\
(3) (\cite[Lemma 4.3]{Ep}) (i) If $M \neq {\Bbb P}$, then $\pi_1(M)$ has no torsion elements.\\
(ii) Suppose $M \neq {\Bbb T}^2$, ${\Bbb K}^2$. If $\alpha, \beta \in \pi_1(M)$ and $\alpha \beta = \beta \alpha$, then $\alpha, \beta \in \langle \gamma \rangle$ for some $\gamma \in \pi_1(M)$. \\
(4) (\cite[p101, lines 5 -10]{Ep}) If $M \neq {\Bbb P}^2$ and every circle component of ${\rm Fr}N$ is essential in $M$, then the inclusion induces a monomorphism $\pi_1(N, x) \to \pi_1(M, x)$ for every $x \in N$.  \\
(5) Suppose $M$ is compact, $X$ is a closed subset of $\partial M$, $X \neq \emptyset$ and $h : M \to M$ is a homeomorphism with $h|_X = id_X$. \\
(i) (\cite[Theorem 3.4]{Ep}) If $M = {\Bbb D}^2$ or ${\Bbb M}$ and $h|_{\partial M} : \partial M \to \partial M$ is orientation preserving, then $h$ is isotopic to $id_M$ rel $X$.  \\
(ii) (\cite[Proof of Theorem 6.3]{Ep}) If $M \neq {\Bbb D}^2$ and $h$ satisfies the following condition ($\ast$), then $h$ is isotopic to $id_M$ rel $X$: \\
($\ast$) $h\ell \simeq \ell$ rel end points for every proper arc $\ell : [0, 1] \to M$ with $\ell(0)$, $\ell(1) \in X$ (we allow that $\ell(0) = \ell(1)$ when $X$ is a single point). 
\end{fact}
\vskip 2mm 
\noindent {\it Comments.}
(4) Consider the universal covering $\pi: \tilde{M} \to M$. By Fact 3.1.(3-i) $\pi^{-1}({\rm Fr}_M N)$ is a union of real lines, half rays and proper arcs. If $M = {\Bbb P}^2$, then $N = {\Bbb P}^2$. \\
(5-ii) $M$ is a disk with $k$ holes, $\ell$ handles (a handle = a torus with a hole) and $m$ M\"obius bands. The assertion is easily verified by the induction on $n = k + \ell + m$, using Fact 3.1.(1) and (5-i), together with the following remarks: \\
(a) When $\# X \geq 2$, we have $h \ell \simeq \ell$ rel. end points even if $\ell(0) = \ell(1)$. \\
(b) If $h$ is (ambient) isotopic to $h_1$ rel. $X$, then $h_1$ also satisfies the condition ($\ast$). \\
(c) Since $M \neq {\Bbb D}^2$, from the condition ($\ast$) it follows that for every component $C$ of $\partial M$, we have $h(C) = C$ and $h$ preserves the orientation of $C$. \\
(d) Let $C_1, \cdots, C_p$ be the components of $\partial M$ which meet $X$. Then $h$ is isotopic rel. $X$ to $h_1$ such that $h_1 = id$ on each $C_i$. Furthermore, $h_1$ satisfies ($\ast$) for $\cup_i C_i$. 

\vskip 3mm 

We also need the following remarks.

\begin{fact}
Suppose $M$ is a connected $2$-manifold and $C$ is a circle component of $\partial M$. If either {\rm (i)} $C \neq \partial M$ or {\rm (ii)} $M$ is noncompact, then $C$ is a retract of $M$. 
\end{fact}

\noindent {\it Comments.}
(ii) Take a half lay $\ell$ connecting $C$ and $\infty$, and consider the regular neighborhood $N$ of $C \cup \ell$. Since $\partial N$ is a real line we can retract $M$ onto $N$ and then onto $C$. \\ 

\begin{fact}
Suppose $M$ is a compact $2$-manifold, $\{ M_i \}$ is a finite collection of compact connected $2$-manifolds such that $M = \cup_i M_i$ and ${\rm Int}M_i \cap {\rm Int}M_j = \emptyset \ (i \neq j)$. \\ {\rm (i)} If $M$ is a disk, then some $M_i$ is a disk. \\ {\rm (ii)} If $M$ is a M\"obius band, then some $M_i$ is a disk or a M\"obius band. \\
{\rm (iii)} If $M$ is an annulus, then some $M_i$ is a disk or an essential annulus in $M$. (If $N$ is a disk with $r$ holes in $M$ ($r \geq 2$), then there exists a disk $D \subset {\rm Int} M$ such that $D \cap N = \partial D \subset \partial N$.)
\end{fact}

\begin{lemma} Suppose $M \neq {\Bbb K}^2$, $C$ is a simple closed curve in $M$ which does not bound a disk or a M\"obius band in $M$, $x \in C$ and $\alpha \in \pi_1(M, x)$ is represented by $C$. If $\beta \in \pi_1(M, x)$ and $\beta^k = \alpha^\ell$ for some $k$, $\ell \in {\Bbb Z} \setminus \{0 \}$, then $\beta \in \langle \alpha \rangle$.
\end{lemma}

\begin{proof}
Attaching $\partial M \times [0, 1)$ to $\partial M \subset M$, we may assume that $\partial M = \emptyset$. Take a covering $p : (\tilde{M}, \tilde{x}) \to (M, x)$ such that $p_{\ast}\pi_1(\tilde{M}, \tilde{x}) = \langle \alpha, \beta \rangle \subset \pi_1(M, x)$. 

If $\tilde{M}$ is noncompact, then by \cite[Lemma 2.2]{Ep} there exists a compact connected 2-submanifold $N$ of $\tilde{M}$ such that $\tilde{x} \in N$ and the inclusion induces an isomorphism $\pi_1(N, \tilde{x}) \to \pi_1(\tilde{M}, \tilde{x})$. Since $\partial N \neq \emptyset$, it follows that $\pi_1(N, \tilde{x}) \cong \langle \alpha, \beta \rangle$ is a free group, so it is an infinite cyclic group $\langle \gamma \rangle$. By Fact.3.1.(2) $\gamma = \alpha^{\pm 1}$, so $\beta \in \langle \alpha \rangle$.

Suppose $\tilde{M}$ is compact. Since ${\rm rank}H_1(\tilde{M}) = 0$ or $1$ and $\pi_1(\tilde{M}) \neq 1$, it follows that $\tilde{M} \cong {\Bbb P}^2$ or ${\Bbb K}^2$ and $M$ is closed and nonorientable. If $\tilde{M} \cong {\Bbb K}^2$ so $\chi(\tilde{M}) = 0$, then $\chi(M) = 0$ and $M \cong {\Bbb K}^2$, a contradiction. Therefore, $\tilde{M} \cong {\Bbb P}^2$ and $\chi(\tilde{M}) = 1$, so $\chi(M) = 1$ and $M \cong {\Bbb P}^2$. We have $\pi_1(M) = \langle \alpha \rangle$. 
\end{proof}

Note that if $M = {\Bbb K}^2$ and $\alpha, \beta$ are represented by the center circles of two M\"obius bands, then $\alpha^2 = \beta^2$, but $\beta \not\in \langle \alpha \rangle$.

\begin{lemma}
Suppose $C$ is a circle component of $\partial M$, $x \in C$ and $\alpha \in \pi_1(M, x)$ is represented by $C$. If $M \neq {\Bbb D}^2$, ${\Bbb M}$ or ${\Bbb S}^1 \times [0, 1] \setminus A$ {\rm (}$A$ is a compact subset of ${\Bbb S}^1 \times \{ 1 \}${\rm )}, then there exists a $\gamma \in \pi_1(M, x)$ such that $\gamma \alpha^n \neq \alpha^n \gamma$ for any $n \in {\Bbb Z} \setminus \{ 0 \}$.
\end{lemma}

\begin{proof}
By the claim below we have a $\gamma \in \pi_1(M, x) \setminus \langle \alpha \rangle$. If $\gamma \alpha^n = \alpha^n \gamma$ for some $n \neq 0$, then by Fact.3.1.(3-ii) $\alpha^n, \gamma \in \langle \beta \rangle$ for some $\beta \in \pi_1(M, x)$ and $\alpha^n = \beta^k$ for some $k \in {\Bbb Z}$. Since $\alpha \neq 1$ and $M \neq {\Bbb P}^2$, by Fact.3.1.(3-i) $k \neq 0$. Hence by Lemma 3.1 $\beta \in \langle \alpha \rangle$ so $\gamma \in \langle \alpha \rangle$, a contradiction. 
\end{proof}

\begin{claim}
Suppose $M$ is a connected $2$-manifold, $C$ is a circle component of $\partial M$, $x \in C$ and $\alpha \in \pi_1(M, x)$ is represented by $C$. If $\pi_1(M, x) = \langle \alpha \rangle$, then $M \cong {\Bbb D}^2$ or ${\Bbb S}^1 \times [0, 1] \setminus A$ for some compact subset $A$ of ${\Bbb S}^1 \times \{ 1 \}$.
\end{claim}

\begin{proof}
First we note that $M$ does not contain any handles or M\"obius bands. In fact if $H$ is a handle or a M\"obius band in $M$, then we can easily construct a retraction $r : M \to H$ which maps $C$ homeomorphically onto $\partial H$ (Fact 3.2), and we have the contradiction $\pi_1(H) = \langle r_{\ast}\alpha \rangle$. In particular, if $M$ is compact then $M$ is a disk or an annulus.

Suppose $M$ is noncompact. It follows that $\partial M$ contains no circle components other than $C$. In fact if $C'$ is a circle in $\partial M \setminus C$, then we can join $C$ and $C'$ by a proper arc $A$ in $M$ and by Fact 3.2 we have a retraction $M \to C \cup A \cup C'$, a contradiction. We can write $M = \cup_{i = 1}^\infty N_i$, where $N_i$ is a compact connected 2-submanifold of $M$, $C \subset {\rm Int}_M N_1$, $N_i \subset {\rm Int}_M N_{i+1}$ and each component of $cl(M \setminus N_i)$ is 
noncompact. We will show that each $N = N_i$ is an annulus. This easily implies the conclusion.

Let $C_1, \cdots, C_m$ be the components of $\partial N \setminus C$. By the above remark $C_j \not\subset \partial M$, so $C_j$ meets a component of $cl(M \setminus N)$. Let $N'$ be a submanifold of $N$ obtained by removing an open color of each $C_j$ from $N$. It follows that $N' \cong N$, ${\rm Fr}N'$ is the union of circles $C_j'$ associated with $C_j$'s, each $C_j'$ is contained in some component $L_j$ of $cl(M \setminus N')$, $cl(M \setminus N') = \cup_j L_j$, and each $L_j$ is noncompact. Since $M$ contains no handles or M\"obius bands (so no one point union of two circles), it follows that $L_j \cap L_j' = \emptyset$ $(j \neq j')$ and $L_j \cap N' = C_j'$. By Fact 3.2 $N'$ is a retract of $M$, so $\pi_1(N', x) = \langle \alpha \rangle$. This implies that $N \cong N'$ is an annulus. 
\end{proof}

The next lemma is a key point in the proof of Theorem 3.1. In \cite[Lemma 6.1]{Ep} the condition \lq\lq the loop $h_t(x)$ is null-homotopic in $M$" is achieved by rotating $x$ along $C$. However, this process does not keep the condition \lq\lq isotopic rel $N$".

\begin{lemma}
Suppose $C$ is a circle component of ${\rm Fr}N$ which does not bound a disk or a M\"obius band in $M$, $h : M \to M$ is a homeomorphism with $h|_N = id_N$ and $h_t: M \to M \ (0 \leq t \leq 1)$ is a homotopy with $h_0 = h$, $h_1 = id_M$. If the following conditions are satisfied, then for any $x \in C$ the loop $m = \{ h_t(x) : 0 \leq t \leq 1 \}$ is null-homotopic in $M$ : \\
{\rm (i)} $M \neq {\Bbb T}^2$, ${\Bbb P}^2$, ${\Bbb K}^2$, \\ 
{\rm (ii)} each circle component of ${\rm Fr}N$ is essential in $M$, \\ 
{\rm (iii)} each component of $N \not\cong {\Bbb D}^2$, ${\Bbb M}$, ${\Bbb S}^1 \times [0, 1] \setminus A$ {\rm (}$A \subset {\Bbb S}^1 \times \{ 1 \}$, compact{\rm )}.
\end{lemma}

\begin{proof}
Let $\alpha = \{ \ell \} \in \pi_1(M, x)$ be represented by $C$ and let $\beta = \{ m \} \in \pi_1(M, x)$. The homotopy $h_t \ell$ implies that $\alpha \beta = \beta \alpha$. Since $M \not\cong {\Bbb T}^2, {\Bbb K}^2$, by Fact.3.1.(3-ii) $\langle \alpha, \beta \rangle \subset \langle \delta \rangle$ for some $\delta \in \pi_1(M, x)$. Since $C$ does not bound a disk or a M\"obius band, by Fact 3.1.(2) $\delta = \alpha^{\pm 1}$ so $\beta = \alpha^k$ for some $k \in {\Bbb Z}$. Let $\alpha_1 = \{ \ell \} \in \pi_1(N, x)$. By Lemma 3.2 there exists a $\gamma = \{ n \} \in \pi_1(N, x)$ such that $\gamma \alpha_1^i \neq \alpha_1^i \gamma$ for any $i \in {\Bbb Z} \setminus \{ 0 \}$ (Figure 2). The homotopy $h_tn$ implies that $\gamma \beta = \beta \gamma$ in $\pi_1(M, x)$. Since $\pi_1(N, x) \to \pi_1(M, x)$ is monomorphic by Fact.3.1.(4), $\gamma \alpha_1^k = \alpha_1^k \gamma$ in $\pi_1(N, x)$ so that $k = 0$ and $\beta = 1$ in $\pi_1(M, x)$.
\end{proof}

\begin{figure}[h]
\caption{The loops $\ell$, $m$ and $n$ in Lemma 3.3.}
\[ \fbox{Figure 2} \]
\end{figure}

\begin{lemma}
Suppose $N \neq \emptyset$, $cl(M \setminus N)$ is compact, each component of ${\rm Fr}N$ is a circle, $h : M \to M$ is a homeomorphism such that $h|_N = id_N$ and $h$ is homotopic to $id_M$. If the following conditions are satisfied, then $h$ is isotopic to $id_M$ rel $N:$ \\ 
{\rm (i)} $M \neq {\Bbb T}^2$, ${\Bbb P}^2$, ${\Bbb K}^2$, \\ 
{\rm (ii)} each component $C$ of ${\rm Fr}N$ does not bound a disk or a M\"obius band, \\
{\rm (iii)} each component of $N \not\cong {\Bbb S}^1 \times [0, 1] \setminus A$ {\rm (}$A \subset {\Bbb S}^1 \times \{ 1 \}$, compact{\rm )}. \end{lemma}

If we assume that $h$ is isotopic to $id_M$, then the condition (iii) is weakened to the condition: \\
(iii)$'$ each component of $N \not\cong {\Bbb S}^1 \times [0, 1]$, ${\Bbb S}^1 \times [0, 1)$. 

\begin{proof}
Let $h_t : h \simeq id_M$ be any homotopy and let $L_1, \cdots, L_m$ be the components of $cl(M \setminus N)$. By Lemma 3.3 the loop $h_t(x) \simeq \ast$ in $M$ for any $x \in {\rm Fr}N = \cup_j {\rm Fr}L_j$. We must find an isotopy $h|_{L_j} \simeq id_{L_j}$ rel ${\rm Fr}L_j$. 

Let $f : [0, 1] \to L_j$ be any path with $f(0)$, $f(1) \in {\rm Fr}L_j$. The homotopy $h_tf$ yields a contraction of the loop $hf \cdot h_t(f(1)) \cdot f^{-1} \cdot (h_t(f(0)))^{-1}$ in $M$. Since $h_t(f(0))$, $h_t(f(1)) \simeq \ast$, it follows that $hf \cdot f^{-1} \simeq \ast$ in $M$. Since $\pi_1(L_j) \to \pi_1(M)$ is monomorphic by Fact.3.1.(4), the loop $hf \cdot f^{-1} \simeq \ast$ in $L_j$, and the desired isotopy is obtained by Fact.3.1.(5-ii).
\end{proof}

Figure 3 illustrates an original idea to prove Lemma 3.4 and Theorem 3.1: Consider the loop $m = n_1 f n_2 f^{-1}$ ($f^{-1}$ is the inverse path of $f$). Any isotopy $h_t : id_M \simeq h$ rel $\{ a, b \}$ induces a homotopy $h_t m : m \simeq n_1 (hf) n_2 (hf)^{-1}$ in $M \setminus \{ a, b \}$. Modify the homotopy $h_t m$ to simplify the intersection of the image of $h_t m$ and ${\rm Fr}N$, and  obtain a homotopy $F : {\Bbb S}^1 \times [0,1] \to M \setminus \{ a, b \}$ shown in Figure 3. The homotopies $F|_{A_1}$ in $N_1$ and $F|_{A_2}$ in $N_2 \setminus \{ a, b \}$ imply that $k_i = 0$ ($i = 1, 2$), and the homotopy $F|_{B_1}$ in $L$ implies that $f \simeq hf$ rel end points in $L$ as required. 

\begin{figure}[h]
\caption{}
\[ \fbox{Figure 3} \]
\end{figure}

\begin{Proof of Theorem 3.1}
We can assume that $X$ is a finite set, since there exists a finite subset $Y$ of $X$ such that $(M, N, Y)$ satisfies the conditions (i) - (iii) in Theorem 3.1. Replacing $h_t$ by $h_1^{-1}h_t$, we may assume that $h_1 = id_M$. 

(I) The case where $M$ is compact : Let $N_1, \cdots, N_n$ be the components of $N$ and $L = cl(M \setminus N)$. Let $K_1, \cdots, K_p$ be the components of $L$ which are disks or M\"obius bands and let $L_1, \cdots, L_q$ be the remaining components. For each $j$ we can write \[ \partial L_j = \left( \cup_{i=1}^{k(j)} A^j_i \right) \cup \left( \cup_{i=1}^{\ell(j)} B^j_i \right) \cup \left( \cup_{i=1}^{m(j)} C^j_i \right), \]
where $A^j_i$'s are the circle components of ${\rm Fr}L_j$, $B^j_i$'s are the components of $\partial L_j$ which contain some arc components of ${\rm Fr}L_j$ and $C^j_i$'s are the remaining components of $\partial L_j$. We choose disjoint collars $E^j_i$ of $A^j_i$ and $F^j_i$ of $B^j_i$ in $L_j$ and set $\hat{A}^j_i = \partial E^j_i \setminus A^j_i$, $\hat{B}^j_i = \partial F^j_i \setminus B^j_i$ and
\[ L_j' = cl(L_j \setminus (\cup_{i=1}^{k(j)} E^j_i) \cup (\cup_{i=1}^{\ell(j)} F^j_i)), \ \ N' = N \cup \left( \cup_{k=1}^p K_k \right) \cup \left( \cup_{j=1}^q \left[ \left( \cup_{i=1}^{k(j)} E^j_i \right) \cup \left( \cup_{i=1}^{\ell(j)} F^j_i \right) \right] \right). \] Note that ${\rm Fr} N' = \cup_{j=1}^q \left[ \left( \cup_{i=1}^{k(j)} \hat{A}^j_i \right) \cup \left( \cup_{i=1}^{\ell(j)} \hat{B}^j_i \right) \right] \subset {\rm Int} M$.
Since ${\cal H}_{\partial}({\Bbb D}) \simeq {\cal H}_{\partial}({\Bbb M}) \simeq \ast$ by Fact 3.1.(5-i), we can isotope $h_0$ rel $N$ to an $h' \in {\cal H}_{N'}(M)$. 

By the construction $(M, N', X, h')$ satisfies the following conditions: \\ 
(1) $N'$ is a 2-submanifold of $M$, every component of ${\rm Fr} N'$ is a circle and $X \subset {\rm Int}_M N'$. \\ 
(2) $h'|_{N'} = id_{N'}$ and $h'$ is isotopic to $id_M$ rel $X$. \\ 
(3) Suppose $C$ is a component of ${\rm Fr} N'$. If $C$ bounds a disk $D$ then $\#(D \cap X) \geq 2$, and if $C$ bounds a M\"obius band $E$ then $E \cap X \neq \emptyset$. \\
(4) If $H$ is an annulus component of $N'$ then $H \cap X \neq \emptyset$. 

To see (3) first note that $M$ is the union of compact $2$-manifolds $N_i$'s, $E^j_i$'s, $F^j_i$'s, $K_k$'s and $L'_j$'s, which have disjoint interiors. Suppose $G$ is a compact connected $2$-manifold in $M$ with $\partial G \subset {\rm Fr} N'$. Since $G \subset {\rm Int} M$ and each $F^j_i$ meets $\partial M$, it follows that $G$ is the union of $N_i$'s, $E^j_i$'s, $L'_j$'s and $K_k$'s contained in $G$. Since $E^j_i$ is an annulus and $L'_j \cong L_j$ is not a disk or a M\"obius band, from Fact 3.3 it follows that
(i) if $G$ is a disk, then $G$ contains a disk which is some $N_i$ or $K_k$ with $K_k \subset G \subset {\rm Int}\, M$, so $\#(G \cap X) \geq 2$, 
(ii) if $G$ is a M\"obius band, then $G$ contains a disk or a M\"obius band which is some $N_i$ or $K_k$ with $K_k \subset G \subset {\rm Int}\, M$, so $G \cap X \neq \emptyset$.

As for (4), $H$ is the union of $N_i$'s $E^j_i$'s, $F^j_i$'s and $K_k$'s contained in $H$, and $H$ contains at least one $N_i$, which is a disk with $r$ holes. If $r \leq 1$ then by the assumption $N_i \cap X \neq \emptyset$. If $r \geq 2$ then we can find a disk $D$ in ${\rm Int} H$ such that $D \cap N_i = \partial D \subset \partial N_i$ (Fact 3.3.(iii)). 
Since $D \subset {\rm Int} N' \subset {\rm Int} M$, $D$ is a union of $N_i$'s and $K_k$'s and we can conclude that it coincides with some $K_k (\subset {\rm Int} M)$, which meets $X$. These imply (4).

It remains to show that $h'$ is isotopic to $id_M$ rel $N'$ under the conditions (1) - (4). 
(i) When $X \subset {\rm Int}N'$, we can apply Lemma 3.4 to the triple $(M \setminus X, N' \setminus X, h'|_{M \setminus X})$. To verify the condition (iii) in Lemma 3.4, note that (a) each component of $N' \setminus X$ takes of the form $H \setminus X$ for some component of $H$ of $N'$, and, in particular, (b) if $H \setminus X \cong {\Bbb S}^1 \times [0, 1)$, then $H$ is a disk and $\#H \cap X \geq 2$, a contradiction. Therefore $h'|_{M \setminus X}$ is isotopic to $id_{M \setminus X}$ rel $N' \setminus X$. Extending this isotopy over $M$ by $id_X$ we have the required isotopy $h' \simeq id_M$ rel $N'$.
(ii) In the case where $X \not\subset {\rm Int} N'$, let $C = \partial N' \cap \partial M$ and consider $(\tilde{M} = M \cup_C C \times [0, 1], \tilde{N} = N' \cup_C C \times [0, 1], X, \tilde{h})$, where $\tilde{h}$ is the extension of $h'$ by $id_{C \times [0,1]}$. Then (a) $X \subset {\rm Int} \tilde{N}$ and $(\tilde{M}, \tilde{N}, X, \tilde{h})$ satisfies (1) - (4), and (b) an isotopy of $\tilde{h}$ to $id_{\tilde{M}}$ rel $\tilde{N}$ restricts to an isotopy of $h'$ to $id_M$ rel $N'$. (Alternatively, we can modify the isotopy of $h'$ to $id_M$ rel $X$ to an isotopy rel $X \cup V$, where $V$ is a neighborhood of $X \cap \partial M$ in $M$. We can replace $X$ so that $X \subset {\rm Int} N'$.) This completes the proof of the case (I).

(II) The case where $M$ is noncompact: Choose a compact connected $2$-submanifolds $L_0$ and $L$ of $M$ such that $h_t(N) \subset {\rm Int}_M L_0 \ (0 \leq t \leq 1)$ and $L_0 \subset {\rm Int}_M L$. Let $N_1 = N \cup cl(L \setminus L_0)$. Since $N_1$ is a subpolyhedron of $L$ with respect to some triangulation of $L$ (cf.\,\cite{Ep}), by Corollary 2.1 we have the principal bundle : ${\cal H}(L)_0 \to {\cal E}(N_1, L)_0^{\ast}$. Let $f_t \in {\cal H}(L)_0$, $f_1 = id_L$, be any lift (= extension) of the path $e_t \in {\cal E}(N_1, L)_0^\ast$ defined by $e_t|_N = h_t|_N$ and $e_t = id$ on $cl(L \setminus L_0)$. 

We can apply the case (I) to $(L, N_1, X_1, f_t)$, $X_1= X \cup cl(L \setminus L_0)$. For the condition (iii) in Theorem 3.1, when $E$ is a component of $cl(L \setminus N_1) = cl(L_0 \setminus N)$, 
(a) if $E \cap {\rm Fr}L_0 = \emptyset$, then $E$ is a component of $cl(M \setminus N)$ and 
(b) if $E \cap {\rm Fr}L_0 \neq \emptyset$, then $E$ contains a component of ${\rm Fr}L_0$ and ${\rm Fr}L_0  \subset cl(L \setminus L_0) \subset X_1$ (it also follows that ${\rm Fr}_L E$ is not connected since $E \cap {\rm Fr}N \neq \emptyset$, so if $E$ is a disk or a M\"obius band, then ${\rm Fr}_L E$ is a disjoint union of arcs).

Therefore we have an isotopy $k_t : L \to L$ rel $N_1$ such that $k_0 = f_0$, $k_1 = id_L$. We can extend $f_t$ and $k_t$ to $M$ by $id$. The required isotopy $h_t'$ is defined by $h_t' = k_tf_t^{-1}h_t$. \qed 
\end{Proof of Theorem 3.1}

\begin{Proof of Corollary 3.1}
Let ${\cal G}_1$ denote the unit path-component of a topological group ${\cal G}$. Theorem 3.1 implies ${\cal H}_N(M) \cap {\cal H}_X(M)_1 = {\cal H}_N(M)_1$. When $M$ is compact, from Fact 2.1 it follows that ${\cal H}_K(M)_0 = {\cal H}_K(M)_1$ for any compact subpolyhedron $K$ of $M$. Since $X$ can be replaced by a finite subset $Y$ of $X$ as in the above proof, we have ${\cal H}_N(M) \cap {\cal H}_X(M)_0 \subset {\cal H}_N(M) \cap {\cal H}_Y(M)_0 = {\cal H}_N(M)_0$. The noncompact case follows from the same argument when we will show that ${\cal H}_K(M)_0$ is an ANR (Propositions 4.1,4.2) in the next section. 
\end{Proof of Corollary 3.1}

\section{The homotopy types of the identity components of homeomorphism groups of noncompact 2-manifolds}

In this final section we will prove Theorem 1.1 and Corollary 1.1.  
Below we assume that $M$ is a {\it noncompact} connected PL 2-manifold and $X$ is a compact subpolyhedron of $M$. 
We set $M_0 = X$ and write as $M = \cup_{i=0}^{\infty} \, M_i$, where for each $i \geq 1$ (a) $M_i$ is a nonempty compact connected PL 2-submanifold of $M$ and $M_{i-1} \subset {\rm Int}_M M_i$, (b) for each component $L$ of $cl \, (M \setminus M_i)$, $L$ is noncompact and $L \cap M_{i+1}$ is connected and (c) $M_1 \cap \partial M \neq \emptyset$ if $\partial M \neq \emptyset$. Taking a subsequence, we have the following cases : \\
\hspace*{4pt} (i) each $M_i$ is a disk, (ii) each $M_i$ is an annulus, (iii) each $M_i$ is a M\"obius band, and \\ \hspace*{4pt} (iv) each $M_i$ is not a disk, an annulus or a M\"obius band. \\
In (ii) the inclusion $M_i \subset M_{i+1}$ is essential, otherwise a boundary circle of $M_i$ bounds a disk component in $cl \, (M \setminus M_i)$, and it contradicts the condition (b).

\begin{lemma}
In the cases $(i) - (iii)$ it follows that 
\begin{tabbing}
$(iii)$ \= \hspace{3pt} (b)$_3$ $M \cong {\Bbb S}^1 \times [0,1] \setminus A$, where $A$ is a nonempty 0-dimensional compact subset of ${\Bbb S}^1 \times \{ 0,1 \}$, \kill
$(i)$ \> $(a)$ $\partial M = \emptyset \Longrightarrow M \cong {\Bbb R}^2$, \\ 
\> $(b)$ $\partial M \neq \emptyset \Longrightarrow M \cong {\Bbb D} \setminus A$, where $A$ is a nonempty 0-dimensional compact subset of $\partial {\Bbb D}$, \\
$(ii)$ \> $(a)$ $\partial M = \emptyset \Longrightarrow M \cong {\Bbb S}^1 \times {\Bbb R}^1$, \\
\> $(b)$ $\partial M \neq \emptyset \Longrightarrow$ \\ 
\> \hspace{3pt} $(b)_1$ $M \cong {\Bbb S}^1 \times [0,1)$ \\ 
\> \hspace{3pt} $(b)_2$ $M \cong {\Bbb S}^1 \times [0,1) \setminus A$, where $A$ is a nonempty 0-dimensional compact subset of ${\Bbb S}^1 \times \{ 0 \}$, \\
\> \hspace{3pt} $(b)_3$ $M \cong {\Bbb S}^1 \times [0,1] \setminus A$, where $A$ is a nonempty 0-dimensional compact subset of ${\Bbb S}^1 \times \{ 0,1 \}$, \\
$(iii)$ \> $(a)$ $\partial M = \emptyset \Longrightarrow M \cong {\Bbb P}^2 \setminus 1pt$, \\
\> $(b)$ $\partial M \neq \emptyset \Longrightarrow M \cong {\Bbb M} \setminus A$, where $A$ is a nonempty 0-dimensional compact subset of $\partial {\Bbb M}$.
\end{tabbing}
\end{lemma}

In the case (ii)(b)$_3$ we may further assume that $M_1$ meets both ${\Bbb S}^1 \times \{ 0 \}$ and ${\Bbb S}^1 \times \{ 1 \}$. 

We choose a metric $d$ on $M$ with $d \leq 1$ and metrize ${\cal H}_X(M)$ by the metric $\rho$ defined by
\[ \rho(f, g) = \sum_{i=1}^{\infty} \, \frac{1}{2^i} \sup_{x \in M_i} \, d(f(x), g(x)). \]

We separate the following two cases: 

\renewcommand{\labelenumi}{(\bf\Roman{enumi})}

\begin{enumerate}
\item $(M, X) \cong ({\Bbb R}^2, \emptyset)$, $({\Bbb R}^2, 1pt)$, $({\Bbb S}^1 \times {\Bbb R}^1, \emptyset)$, $({\Bbb S}^1 \times [0, 1), \emptyset)$, $({\Bbb P}^2 \setminus 1pt, \emptyset)$. 
\item $(M, X)$ is not Case (I).
\end{enumerate}
\vskip 2mm
\noindent {\bf Case (II)}: First we treat Case (II) and prove the following statements: 

\begin{proposition}
In Case {\rm (II)}, we have {\rm (1)} ${\cal H}_X(M)_0 \simeq \ast$ and {\rm (2)} ${\cal H}_X(M)_0$ is an ANR.
\end{proposition}

We use the following notation: For each $j \geq 1$ let $U_j = {\rm Int}_M M_j$ and $L_j = {\rm Fr}_M M_j$, and for each $j > i \geq k \geq 0$ let ${\cal H}_{k,j} = {\cal H}_{M_k \cup (M \setminus U_j)}(M)_0$, ${\cal U}^i_{k,j} = {\cal E}_{M_k}(M_i, U_j)^{\ast}_0$ and let $\pi^i_{k,j} : {\cal H}_{k,j} \to {\cal U}^i_{k,j}$ denote the restriction map, $\pi^i_{k,j}(h) = h|_{M_i}$. 

\begin{lemma}
$(1)$ ${\cal H}_{k,j} \cong {\cal H}_{M_k \cup L_j}(M_j)_0$ is an AR . \\
$(2)$ The map $\pi^i_{k,j} : {\cal H}_{k,j} \to {\cal U}^i_{k,j}$ is a principal bundle with the structure group ${\cal H}_{k,j} \cap {\cal H}_{M_i}(M) \cong {\cal G}^i_{k,j} \equiv {\cal H}_{M_k \cup L_j}(M_j)_0 \cap {\cal H}_{M_i \cup L_j}(M_j)$ (under the restriction map).
\\
$(3)$ ${\cal U}^i_{k,j}$ is an open subset of ${\cal E}_{M_k}(M_i, M)^{\ast}_0$, $cl \, {\cal U}^i_{k,j} \subset {\cal U}^i_{k,j+1}$ and ${\cal E}_{M_k}(M_i, M)^{\ast}_0 = \cup_{j > i} \, {\cal U}^i_{k,j}$. 
\end{lemma}

\begin{proof}
The statement (1) follows from Fact 2.1 and Lemma 2.1.(ii), and (2) follows from Corollary 2.1. For (3), note that ${\cal E}_{M_k}(M_i, M)^{\ast}_0$ is path connected (Proposition 2.1) and each $f \in {\cal E}_{M_k}(M_i, M)^{\ast}_0$ is isotopic to the inclusion $M_i \subset M$ in a compact subset of $M$. 
\end{proof}

\begin{lemma}
In Case {\rm (II)}, for each $j > i > k \geq 0$, {\rm (a)} ${\cal G}^i_{k,j}$ is an AR, {\rm (b)} the restriction map $\pi^i_{k,j} : {\cal H}_{k,j} \to {\cal U}^i_{k,j}$ is a trivial bundle and {\rm (c)} ${\cal U}^i_{k,j}$ is also an AR. \end{lemma}

\begin{proof}
Once we show that ($\ast$) ${\cal G}^i_{k,j} = {\cal H}_{M_i \cup L_j}(M_j)_0$, then (a) the fiber ${\cal G}^i_{k,j}$ is an AR by Fact 2.1 and Lemma 2.1.(ii), so (b) the principal bundle has a global section and it is trivial and (c) follows from Lemma 4.2.(1). It remains to prove ($\ast$).  \\
(1) The cases (i)(a), (ii)(a), (iii)(a), (ii)(b)$_1$ and (iv) (under the condition (II)): 
We can apply Theorem 3.1 to $(\tilde{M}_j = M_j \cup_{L_j} L_j \times [0, 1], \tilde{M}_i = M_i \cup_{L_j} L_j \times [0, 1], \tilde{M}_k = M_k \cup_{L_j} L_j \times [0, 1])$. We can verify the conditions (ii) and (iii) in Theorem 3.1 as follows: (ii) By the assumption $(M_i, X) \not\cong ({\Bbb D}, \emptyset)$, $({\Bbb D}, 1pt)$, $({\Bbb S}^1 \times [0,1], \emptyset)$, $({\Bbb M}$, $\emptyset)$ for each $i \geq 1$. (iii) If $H$ is a component of $cl_{\tilde{M}_j}(\tilde{M}_j \setminus \tilde{M}_i) = cl_M(M_j \setminus M_i)$, then $H$ contains a component of $L_j$. (Also, $H$ meets both $M_i$ and $L_j$ (if $H \cap L_j = \emptyset$ then $H$ is a compact component of $cl(M \setminus M_i)$, a contradiction.), so ${\rm Fr}_{\tilde{M}_j} H$ is not connected. Hence if $H$ is a disk or a M\"obius band, then ${\rm Fr}_{\tilde{M}_j} H$ is a disjoint union of arcs.) By Corollary 3.1 (Compact case) it follows that ${\cal H}_{\Tilde{M}_k}(\tilde{M}_j)_0 \cap {\cal H}_{\tilde{M}_i}(\tilde{M}_j) = {\cal H}_{\tilde{M}_i}(\tilde{M}_j)_0$ and this implies ($\ast$).

(2) The cases (i)(b), (ii)(b)$_{3}$ and (iii)(b): 
Since $M_1 \cap \partial M \neq \emptyset$ and $M_1$ meets both ${\Bbb S}^1 \times \{ 0 \}$ and ${\Bbb S}^1 \times \{ 1 \}$ in the case (ii)(b)$_3$, it follows that $cl(M_j \setminus M_i)$ is a disjoint union of disks, thus ${\cal H}_{M_i \cup L_j}(M_j) = {\cal H}_{M_i \cup L_j}(M_j)_0$ by Fact 3.1.(5-i). This implies ($\ast$). 

(3) The remaining case (ii)(b)$_2$: 
It follows that (a) $cl(M_j \setminus M_i)$ is a disjoint union of disks $D_k$ and an annulus $H$ and (b) $D_k \cap M_i$ is an arc, $D_k \cap L_j$ is a disjoint union of arcs ($\neq \emptyset$) and $H \cap M_i$, $H \cap L_j$ are the boundary circles of $H$. Since $N = M_i \cup (\cup_k D_k)$ is an annulus and $N \cap L_j \neq \emptyset$, from Theorem 3.1 it follows that ${\cal H}_{M_k \cup L_j}(M_j)_0 \cap {\cal H}_{N \cup L_j}(M_j) = {\cal H}_{N \cup L_j}(M_j)_0$. Each $f \in {\cal G}^i_{k,j}$ is isotopic rel $M_i \cup L_j$ to $f' \in {\cal H}_{N \cup L_j}(M_j)$. Since $f'$ is isotopic to id rel $M_k \cup L_j$, it follows that $f' \in {\cal H}_{N \cup L_j}(M_j)_0$ and so $f \in {\cal H}_{M_i \cup L_j}(M_j)_0$. This completes the proof.
\end{proof}

\begin{lemma}
In Case {\rm (II)}, for each $i \geq k \geq 0$, \\
{\rm (a)} ${\cal E}_{M_k}(M_i, M)^{\ast}_0$ is an AR, \\ 
{\rm (b)} the restriction map 
$\pi : {\cal H}_{M_k}(M)_0 \longrightarrow {\cal E}_{M_k}(M_i, M)^{\ast}_0$ 
is a trivial principal bundle with fiber ${\cal H}_{M_i}(M)_0$, \\
{\rm (c)} ${\cal H}_{M_k}(M)_0$ strongly deformation retracts onto ${\cal H}_{M_i}(M)_0$. 
\end{lemma}

\begin{proof}
By Lemma 4.3.(c) each ${\cal U}^i_{k,j} \ (j > i)$ is an AR. Thus by Fact 4.2.(3) ${\cal E}_{M_k}(M_i, M)^{\ast}_0$ is also an AR and it strongly deformation retracts onto the single point set $\{ M_i \subset M \}$. Hence the principal bundle 
\[ {\cal G}_k^i \equiv {\cal H}_{M_k}(M)_0 \cap {\cal H}_{M_i}(M) \subset {\cal H}_{M_k}(M)_0 \longrightarrow {\cal E}_{M_k}(M_i, M)^{\ast}_0 \] 
is trivial and ${\cal H}_{M_k}(M)_0$ strongly deformation retracts onto the fiber ${\cal G}_k^i$. In particular, ${\cal G}_k^i$ is connected and ${\cal G}_k^i = {\cal H}_{M_i}(M)_0$.
\end{proof}

\begin{Proof of Proposition 4.1.(1)} 
By Lemma 4.4.(c), for each $i \geq 0$ there exists a strong deformation retraction $h^i_t$ ($0 \leq t \leq 1$) of ${\cal H}_{M_i}(M)_0$ onto ${\cal H}_{M_{i+1}}(M)_0$. A strong deformation retraction $h_t$ ($0 \leq t \leq \infty$) of ${\cal H}_X(M)_0$ onto $\{ id_M \}$ is defined as follows:
\begin{eqnarray*}
 h_t(f) \, &=& h_{t-i}^ih_1^{i-1} \cdots h_1^0(f) \ \ \ (f \in {\cal H}_X(M)_0, \ i \geq 0, \ i \leq t \leq i+1) \\
 h_\infty(f) &=& id_M.
\end{eqnarray*}
Since ${\rm diam} \, {\cal H}_{M_i}(M)_0 \leq 1/2^i \to 0$, the map $h : {\cal H}_X(M)_0 \times [0, \infty] \to {\cal H}_X(M)_0$ is continuous. 

(In the cases (i), (ii) and (iii), the same conclusion follows from Lemma 2.1.(i)(ii) by taking the end compactification of $M$.) \qed
\end{Proof of Proposition 4.1.(1)} 

For the proof of Proposition 4.1.(2), we will apply Hanner's criterion of ANRs: 

\begin{fact} (\cite{Han}) 
A metric space $X$ is an ANR iff for any $\varepsilon > 0$ there is an ANR $Y$ and maps $f : X \to Y$ and $g : Y \to X$ such that $gf$ is $\varepsilon$-homotopic to $id_X$. 
\end{fact}

\begin{Proof of Proposition 4.1.(2)} 
By Lemma 4.4.(b) and Proposition 4.1.(1) for each $i \geq 1$ we have the trivial principal bundle
\[ {\cal H}_{M_i}(M)_0 \subset {\cal H}_X(M)_0 \stackrel{\pi}{\longrightarrow} {\cal E}_X(M_i, M)^{\ast}_0 \ \ {\rm with} \ \ {\cal H}_{M_i}(M)_0 \simeq \ast. \]
It follows that $\pi$ admits a section $s$, and the map $s \pi$ is fiber preserving homotopic to $id_{{\cal H}_X(M)_0}$ over ${\cal E}_X(M_i, M)^{\ast}_0$. Since each fiber of $\pi$ has ${\rm diam} \leq 1/2^i$, this homotopy is a $1/2^i$-homotopy. Since ${\cal E}_X(M_i, M)^{\ast}_0$ is an ANR (Proposition 2.1), by Fact 4.1 ${\cal H}_X(M)_0$ is also an ANR. \qed
\end{Proof of Proposition 4.1.(2)} 

\noindent {\bf Case (I)}: The next statements follow from Lemma 2.1.(iii) and Fact 2.1 by taking the end compactification of $M$. 

\begin{proposition}
In Case {\rm (I)}, we have {\rm (1)} ${\cal H}_X(M)_0 \simeq {\Bbb S}^1$ and {\rm (2)} ${\cal H}_X(M)_0$ is an ANR.
\end{proposition}

Theorem 1.1 follows from Propositions 4.1, 4.2, and Corollary 1.1 now follows from the following characterization of $\ell_2$-manifold topological groups.

\begin{fact} (\cite{DT})
A topological group is an $\ell_2$-manifold iff it is a separable, non locally compact, completely metrizable ANR. 
\end{fact}

\begin{Proof of Corollary 1.1} Since $M$ is locally compact and locally connected, ${\mathcal H}(M)$ is a topological group and ${\mathcal H}_X(M)$ is a closed subgroup of ${\mathcal H}(M)$. Since $M$ is locally compact and second countable, ${\mathcal H}(M)$ is also second countable. A complete metric $\rho$ on ${\cal H}(M)$ is defined by 
\[ \rho(f, g) = d_{\infty}(f, g) + d_{\infty}(f^{-1}, g^{-1}), \hspace{0.5cm} d_{\infty}(f, g) = \sum_{n = 1}^{\infty} \frac{1}{2^n} \sup_{x \in M_n} d(f(x), g(x)) \]
for $f, g \in {\mathcal H}(M)$, where $d$ is a complete metric on $M$ with $d \leq 1$. Since ${\mathcal H}_X(M)_0 \cong {\mathcal H}_X(M)_0 \times s$ \cite{Ge}, ${\mathcal H}_X(M)_0$ is not locally compact. Finally, by Propositions 4.1, 4.2 ${\mathcal H}_X(M)_0$ is an ANR. This completes the proof. \qed
\end{Proof of Corollary 1.1}

\end{document}